\RequirePackage{fix-cm}
\documentclass{article}

\usepackage{amsfonts}
\usepackage{amsmath,amssymb}
\usepackage{graphicx}
\usepackage{epstopdf}
\usepackage{float} 
\usepackage{xfrac} 
\usepackage{array}
\usepackage{hyperref}
\usepackage{cleveref}
\usepackage[format=plain,labelfont=bf,
            textfont=it,font=small]{caption}
\usepackage{geometry}

\newcommand{\email}[1]{\texttt{#1}}

\newcommand{\Xag}{\mathbf{X}}
\newcommand{\Rmat}{\mathbb{R}}
\newcommand{\SPS}{Single}
\newcommand{\MPA}{Mixed1}
\newcommand{\MPB}{Mixed2}
\newcommand{\DPS}{Double}
\newcommand{\Kur}{Benchmark 2}
\newcommand{\LCO}{Benchmark 1}
\newcommand{\CC}{Benchmark 3}
\newcommand{\ie}{\textit{i.e.}}

\newcommand{\figcom}[2]{\ifx&#1&\ifx&#2&\else({\it #2})\fi\else\ifx&#2&(#1)\else(#1, {\it #2})\fi\fi}

\title{Mixed-Precision in adaptive Runge-Kutta method for large ODE systems\footnote{This manuscript has been submitted to the SIAM Journal on Scientific Computing and is currently under revision.}}

\author{Mouhamad Al-Sayed Ali\thanks{URMAR Univ. Rennes, Rennes, 35000, France. (\email{mouhamad.alsayedali@univ-rennes.fr}).}
\and Samuel Bernard\thanks{MUSICS Team, Inria Lyon, Villeurbanne, France and CNRS, Université Claude Bernard Lyon 1, ICJ UMR5208, Inria.
  (\email{samuel.bernard@inria.fr}).}
\and Arsène Marzorati\thanks{BioTIC and MUSICS Teams, Inria Lyon, Villeurbanne, France and INSA de Lyon, Université Claude Bernard Lyon 1, ICJ UMR5208, Inria. (\email{arsene.marzorati@inria.fr}).}
\and Jonathan Rouzaud-Cornabas\thanks{BioTIC Team, Inria Lyon, Villeurbanne, France and INSA Lyon, Inria, CITI, UR3720, 69621 Villeurbanne, France(\email{jonathan.rouzaud-cornabas@inria.fr}).}
}

\begin{document}
\maketitle 
\begin{abstract}
		
		Mixed-precision methods combine low and high precision arithmetics to exploit low precision computational speed and high precision accuracy. Large ODE systems that contain many heterogeneous interactions lead to a high computational cost that could be tackled with mixed-precision solvers. 
		We tested mixed-precision versions of the Bogacki-Shampine 3(2) Runge-Kutta pair over three benchmark systems: coupled linear oscillators, the Kuramoto model and a circadian clock model. Our study is performed in a way that can be adapted to any finite-precision format, software architecture and numerical scheme. 
		We found that mixed-precision solvers can preserve most of the high-precision solver accuracy under a wide range of solver tolerances. Moreover, mixed-precision solver accuracy improves with system size, reaching levels equivalent to high-precision solvers in small system size. We also observed that mixed-precision arithmetic does not impact the number of evaluation in a way that balances the benefit of fast operations in low precision. 
		Taken together, these results show that mixed-precision methods can offer significant computational speed-up at little or no loss of accuracy in large coupled ODE systems.	
\end{abstract}

\section{Introduction}\label{S:Introduction}
	
	Numerical methods for dynamical systems such as ordinary differential equations (ODEs) are widely used in many scientific fields such as biology \cite{MURRAY-2002}, physics or engineering \cite{HIRSCH-2013,STROGATZ-2018} and for solving discretized partial  differential equations \cite{BAKER-2000}. This is why this kind of problem has been studied for a long time and numerous numerical tools for solving it have been developed \cite{BUTCHER-2016,HAIRER-1993}. 
	However, the applicability of these systems is often limited by their complexity, \ie{} their size and connectivity. For example, in computational biology with bacteria ecosystems \cite{NAGRAJAN-2022}, in meteorology \cite{KURTH-2018,PALMER-2014} or more generally in geophysics models \cite{ACKMANN-2022}, large systems have to be solved, requiring large amount of data to be stored and manipulated. 
	
	Computational needs of large ODE systems have led to the development and specification of optimized and well-tuned solvers \cite{CAVAGLIERI-2015,ICHIMURA-2018,KENNEDY-2000}. Dimension reduction techniques have been developed, such as reduced order modeling \cite{ANTOULAS-2004,DAR-2023,LUCIA-2004} which involves transforming the original large system into a much smaller one, or mean-field theory \cite{AOKI-2014,BORDENAVE-2010,BORTOLUSSI-2016} which averages the effective interactions. However, while these methods are effective, they require additional assumptions, for example homogeneity, sparsity or locality in interactions. These hypotheses may be restrictive in some real cases, such as in the study of biological networks \cite{PAUL-2020-BIONET}. Our study addresses a more general case including fully pairwise interactions.
	To avoid unnecessary, constraints we explore the possibility of using mixed-precision methods for solving large coupled ODE systems. Mixed-precision methods consist in using several arithmetic precision formats within a single algorithm in order to increase computation speed while reducing the space required to store the model. Indeed, there are many formats for encoding numbers on computers and performing numerical computations. The most widespread formats for real numbers are the floating-point standards of the IEEE 754 \cite{ZURAS-2008} with half (16 bits), single (32 bits)  and double (64 bits) precision. These finite representations influence arithmetic, accuracy and computing speed. Accuracy depends on the {\it machine epsilon} which is the smallest value summable to 1 and speed (storage, memory displacement and arithmetic operations) is directly impacted by the number of bits used for the encoding.
	
    Recent studies have demonstrated the advantages of mixed-precision methods across various scientific fields. These approaches combine different numerical precision formats to optimize both accuracy and efficiency \cite{LEGRAND-2013-SPFP}: higher-precision formats are selectively applied to the critical components of algorithms to preserve or enhance accuracy \cite{CHIANG-2017-RIGFP}, while lower-precision formats help reduce computational and storage costs due to their smaller bit sizes \cite{HAYFORD-2024}.

	Mixed-precision approaches have been particularly successful in linear algebra, where algorithms such as GMRES and LU factorization have been effectively re-implemented \cite{HIGHAM-2022}. They are also widely used in deep learning, enhancing both training \cite{MICIKEVICIUS-2017-MPTRAIN} and inference phases through quantization of weights and activations \cite{RAKKA-2022-MPNN}. Beyond these areas, mixed-precision computing has found applications in computational fluid dynamics \cite{BROGI-2024-MPCFD}, climate forecasting \cite{SAFFIN-2020-REDUCED}, and other scientific disciplines \cite{KASHI-2024-MPSURVEY}.
	
	Of the many numerical schemes available \cite{HAIRER-1993}, few studies have proposed to adapt the algorithms with mixed-precision. In \cite{BURNETT-2022, GRANT-2020}, \textit{Grant and al.} develop a complete framework for implicit additive Runge-Kutta (RK) methods. While \textit{Croci and al.} make a complete study of the accuracy for mixed-precision versions of stabilized explicit Runge-Kutta-Chebyshev methods in \cite{CROCI-2022}.
	
	In this work, we study the impact of mixed-precision on the accuracy when solving nonlinear ODE systems. Contrary to previous work, we choose another approach working with an adaptive explicit scheme. The stability criterion for explicit schemes is closely linked to the size of the discretization time step. This parameter can be set before computing the solution, or it can be adjusted at runtime to control the stability of the numerical scheme and limit the number of operations.
	The Bogacki-Shampine 3(2) pair \cite{BOGACKI-1989} was designed to optimize the intersection of the stability regions of its two low order RK schemes and limit the computational cost by using embedded schemes and the First Same As Last ({\it FSAL}) property. For its simplicity, efficiency and economical aspects, the study is conducted on this numerical scheme.

	We first present the general form of ODEs system which could take advantage of mixed-precision methods. Next, three benchmarks satisfying this framework are introduced: (1) linear oscillators with a linear coupling term, (2) the Kuramoto system \cite{KURAMOTO-1984}, a standard representation of synchronization with a non-linear interaction term, and (3) a circadian clock model (\Cref{SbS:ODE-Systems}). The latter combines a version of the Goodwin model \cite{Goodwin-1965} describing a negative feedback oscillators in cellular system, such as circadian rhythms, and a version of the FitzHugh-Nagumo \cite{FITZHUGH-1961,NAGUMO-1962} model which was originally developed to describe the electro-physiology of neurons. We use it here to describe the balance between different protein complexes driving the sharp transition to mitosis during the cell cycle. 
	We then briefly recall the structure of the Bogacki-Shampine 3(2) pair and detail the introduction of mixed-precision into the scheme (\Cref{SbS:MP-Numerical-Scheme}). We set how low precision is distributed inside the system evaluation and among the stages of the numerical scheme. We implement in Matlab four versions of the numerical scheme, one in single, one in double precision and two in mixed-precision (\Cref{SbS:MP-RK}). And finally we define the metrics used to exploit the results of the different solvers on the benchmarks (\Cref{SbS:Error-Metric}). In particular, we introduce a performance proxy as run-times on Matlab are not representative of actual performance. Indeed, computations in Matlab are enforced to satisfy the required arithmetic precision, but this entails a loss of global performance (\Cref{SbS:Perf-Criteria}).
	We organise the results into four main arguments. First, we focus on how to include mixed-precision into the numerical scheme. The distribution and number of operations performed in low precision leverage performance. They necessarily modify the potential acceleration of a faster computation, but they also affect the number of steps computed by the solver (\Cref{SbS:Perf-Proxy}).
	Second, we look at accuracy with respect to mixed-precision and system sizes. Mixed-precision solvers appear to be more robust to tolerance constraints than the single precision one and, in some cases, can reach an accuracy similar to that of the double precision solver (\Cref{SbS:Accuracy-Tolerance}). In addition, all solvers benefit from the increase in system size. Finally, we show that the error control fails for the single precision solver at smaller tolerances, whereas the mixed-precision solvers better match the expected theoretical behaviour even at strict tolerances (\Cref{SbS:Local-Error}).

\section{Methods}\label{S:Methods}
	
	In this section, we present the general framework and three ODE systems which are used as benchmarks. Then, we introduce the numerical scheme applied on the systems, and we explain how we insert the mixed-precision in its calculations. We end this section with a presentation of the metrics and the performance proxy.
	
	\subsection{ODE system: General model and benchmarks}\label{SbS:ODE-Systems}
	\paragraph{General model} We consider a system describing a population of $N$ heterogeneous agents, with $N \gg 1$, described by the variable $X \in \Rmat^{d N}$. Each agent is identified with the variable $\Xag_i$ in $\Rmat^{d}$, considering the following extraction $\Xag_i=\Big(X_{(i-1)d+1},...,X_{id}\Big)$. The dynamical dependence  is split into two parts, an agent-centered term $F_i:\Rmat\times\Rmat^d\to\Rmat^d$ and a term $G_{ij}:\Rmat\times\Rmat^d\times\Rmat^d\to\Rmat^d$ accounting for complex pairwise interactions, which are weighted by a vector $M_{ij} \in \Rmat^d$. The state of each agent follows the differential equation:
	\begin{equation}
		\dot{\Xag_i}(t) = F_i(t,\Xag_i) + \sum_{j=1}^N M_{ij}(t,\Xag_i)\odot G_{ij}(t,\Xag_i,\Xag_j), \,i\in{1,...,N},
	\end{equation} 
	where $\odot$ is the Hadamard product (or element-wise product).
	There is no additional assumption on the form of $G_{ij}$ except that all the $N$ agents can fully interact heterogeneously, leading to an evaluation of $N^2$ terms in the right-hand side.

	We build the agent-centered term $F: \Rmat \times \Rmat^{dN} \to \Rmat^{dN}$ for the whole population by extending the functions $F_i$ blockwise. In the same way, the interaction terms $G_{ij}$ and $M_{ij}$ extend to a population interaction function $G: \Rmat \times \Rmat^{dN}  \to \Rmat^{dN\times N}$ and an array of vectors of weights $M \in \Rmat^{dN\times N}$.  We denote as a convolution product the weighting of the interactions ($[M*G]: \Rmat \times \Rmat^{dN} \to \Rmat^{dN}$).  The system of equations describing the general framework for the whole population is 
	\begin{equation}\label{EQ:General-Model}
		\dot{X}(t) = F(t,X) + [M*G](t,X).
	\end{equation}

	Now, we present three benchmarks conforming to \Cref{EQ:General-Model}. For each benchmark, we explicitly define the functions $F_i$ and $G_{ij}$ and the weight vector $M_{ij}$. All the parameter values are listed in  \Cref{Tab:Parameters-Values-LCO,Tab:Parameters-Values-Kur,Tab:Parameters-Values-CC} of \Cref{AS:Parameters}. 
	
	\paragraph{Coupled linear oscillators (\LCO{})} \LCO{} is a system of harmonic oscillators coupling by a linear term. 
	\begin{equation}\label{EQ:LCO}
		\begin{aligned} 
			F^{LCO}_i:&\; \mathbb{R}^2 \to \mathbb{R}^2 & G^{LCO}_{ij}:&\; \mathbb{R}^2 \times \mathbb{R}^2 \to \mathbb{R}^2 \\
			&\Xag_i \mapsto \begin{pmatrix}\Xag_{i,2}\\ -\Xag_{i,1}\end{pmatrix}& &(\Xag_i,\Xag_j)\mapsto \begin{pmatrix}\Xag_{j,1}-\Xag_{i,1}\\ 0\end{pmatrix}\\
		\end{aligned}
	\end{equation}
	\[M_{ij}=(\frac{1}{N},0)^T \]
	
	System \Cref{EQ:LCO} posses an analytic solution.
	
	\paragraph{Kuramoto model (\Kur{})} \Kur{} is the classic Kuramoto model \cite{KURAMOTO-1984}.
	\begin{equation}\label{EQ:Kuramoto}
		\begin{aligned} 
			F^{Kur}_i:&\; \mathbb{R} \to \mathbb{R} & G^{Kur}_{ij}:&\; \mathbb{R} \times \mathbb{R} \to \mathbb{R} \\
			&\Xag_i \mapsto \omega_i& &(\Xag_i,\Xag_j)\mapsto K\sin(\Xag_j-\Xag_i)\\
		\end{aligned}
	\end{equation}
	\[M_{ij}=\frac{1}{N}\]
	Parameters are:  $\omega_i$, the natural frequency of the $i$th oscillator, and $K$, the coupling constant.
	The $\omega_i$ are generated with a centered normal distribution with standard deviation $\sigma$. $\sigma$ and $K$ are changed for each test respectively between $[0,1]$ and $[0, 3\sigma]$ (\Cref{ASbS:Parameters-Kur}).
	
	\paragraph{Circadian clock/cell cycle model (\CC{})} \CC{} combines two nonlinear oscillators coupled through a nonlinear term.
	
	\begin{equation}\label{EQ:CC}
		\begin{aligned} 
			F^{CC}_i:&\; \mathbb{R}^4 \to \mathbb{R}^4  \\
			&\Xag_i \mapsto 
			\Bigg(
			 \frac{k_0 \theta^h }{\theta^h+\Xag_{i,2}^h}\bigg(a\Xag_{i,1}^2+1\bigg) - k_1 \Xag_{i,1}, \, 
			k_2 (\Xag_{i,1} - \Xag_{i,2}), \, \\
			&\quad\quad\quad\quad \Xag_{i,3} \bigg( 1 - \frac{\Xag_{i,3}^2}{3} \bigg) - \Xag_{i,4} 
			+ \mathrm{I}_0 \bigg( 1 - \frac{k_3^2}{k_3^2 + \Xag_{i,1}^2} \bigg), \,
			\epsilon \big( \Xag_{i,3} + b - c \Xag_{i,4} \big)
			\Bigg)^{\mathrm{T}} \\
			G^{CC}_{ij}:&\; \mathbb{R}^4 \times \mathbb{R}^4 \to \mathbb{R}^4\\
			&(\Xag_i, \Xag_j) \mapsto 
			\begin{pmatrix}
				\frac{k_0 \theta^h a}{\theta^h+\Xag_{i,2}^h} K \arctan(\Xag_{j,1} - \Xag_{i,1}), \, 0, \, 0, \, 0
			\end{pmatrix}^{\mathrm{T}}\\
		\end{aligned}
	\end{equation}
	\[M_{ij}=(\frac{1}{N},0,0,0)^T \]
	
	Parameter values are fixed for $k_0=2$, $k_2=0.144832$, $k_3=2$, $a=2$, $b=0.7$, $c=0.8$, $h=4$, $\epsilon = 0.228249$.
	Parameter $k1$ takes random value following a normal distribution with mean $\Bar{k}_1=0.339278$ and standard deviation $\sigma_1=0.090909$. 
	The parameter $\theta^h = \frac{k_1}{k_0-k_1}$. The coupling coefficient $K$ and the stimulus $I_0$ are changed for each test (\Cref{ASbS:Parameters-CC}).

	\subsection{Numerical scheme calculations and mixed-precision}\label{SbS:MP-Numerical-Scheme}
	
	\paragraph{The Bogacki-Shampine 3(2) pair} We chose an adaptive scheme with embedded Runge-Kutta methods \cite{BOGACKI-1989}.
	It combines a third order method, with three stages, and a second order one with four stages \Cref{Tab:Butcher-O23}. At the {\it n-th} step of the scheme we denote by $X^n$ the third order solution and $\Tilde{X}^n$ the second order one. To get the solutions at the next step, the following four stages ($K_{\ell}$) are computed:
	\begin{equation}\label{EQ:Stage-ODE23}
		K_{\ell} = F\big(t+c_{\ell}h,X^n+h\sum_{m=1}^{\ell-1}a_{\ell m}K_m\big)+[M*G]\big(t+c_{\ell}h,X^n+h\sum_{m=1}^{\ell-1}a_{\ell m}K_m\big), \,\ell=1,2,3,4.
	\end{equation}
	
	The solutions at the next step $t+h$ are given by:
	\[X^{n+1} = X^n + h\sum_{\ell=1}^3b_{\ell}K_{\ell}.\]
	\[\Tilde{X}^{n+1} = X^n + h\sum_{\ell=1}^4\Tilde{b}_{\ell}K_{\ell}.\]
	
	\begin{table}[htbp]
		\centering
		$\begin{array}
			{c|cccc}
			c_{\ell}&a_{\ell 1} & a_{\ell 2} & a_{\ell 3}\\
			\hline
			0\\
			\sfrac{1}{2} & \sfrac{1}{2}\\
			\sfrac{3}{4} &0 &\sfrac{3}{4} \\
			1& \sfrac{2}{9} & \sfrac{1}{3} & \sfrac{4}{9} \\
			\hline
			b_m & \sfrac{2}{9} & \sfrac{1}{3} & \sfrac{4}{9} & 0 \\
			\Tilde{b}_m & \sfrac{7}{24} & \sfrac{1}{4} & \sfrac{1}{3} & \sfrac{1}{8} \\
		\end{array}$
		\caption{Butcher table for the Bogacki-Shampine 3(2) pair \cite{BOGACKI-1989}.}
		\label{Tab:Butcher-O23}
	\end{table}
	The fourth stage, used only in the second order solution, corresponds to the evaluation of the first stage in the next step. Thus, if the step is accepted, the first stage is already computed (First Same As Last method). This method minimizes the computational cost by two aspects: the embedded stages of the schemes and the FSAL property. The latter applies only when a step is accepted, otherwise the first stage of the next step (stage $K_4$) must be recomputed while the current step is rejected.
	
	\paragraph{Local relative error estimation} The specificity of adaptive solvers is the modification of the step size ($h$) to control the local error, \ie{} the error introduced in a single step. The step is accepted only if the approximation of the local error meets a tolerance criterion. This approximation is computed as the difference between the two numerical solutions $X^{n+1}$ and $\Tilde{X}^{n+1}$.
	The tolerance criteria are defined with the absolute ($Ab$) and relative ($Rel$) tolerances.
	
	The relative local error approximation is done with a normalization vector of weights $W$ which corresponds to the norm of the higher order solution vector and enables to switch from relative to absolute tolerance criterion if this norm of the solution is too small. Then, using the stages (see \Cref{EQ:Stage-ODE23}), we compute the difference vector $V_{BS}$ between the two solutions of the embedded schemes. This difference is weighted with $W$. Finally, the local relative error, $E_{BS}^n$, is approximated with the maximum of the absolute value of $V_{BS}$ elements:
	\begin{equation}\label{EQ:ERROR-O23}
			\begin{aligned}
		W_{k}=&\max(|X_{k}^n|, |X_{k}^{n+1}|,\tfrac{Ab}{Rel}),& k \in \{1,...,dN\}.\\
		V_{BS,k}=&\frac{|X^{n+1}_{k}-\Tilde{X}^{n+1}_{k}|}{W_k},&  k \in \{1,...,dN\}.\\
		E_{BS}^n =& \max_{k}(V_{BS,k}).
	\end{aligned}
	\end{equation}
	We chose the metrics for \Cref{EQ:ERROR-O23} based on published specifications \cite{SHAMPINE-1997}, but others choices are possible \cite[Section II.4]{ HAIRER-1993}.
	
	The step is accepted if and only if the local relative error ($E_{BS}^n$) is smaller than the relative tolerance. The error approximation is also used within the correction of the step size. 
	We recall that a numerical scheme of order $p$ has a local truncation error of magnitude $O(h^{p+1})$. Here, if a step is accepted, the solution of the third order ($p=3$), $X^{n+1}$ is taken as the solution of the numerical scheme.
	
	\paragraph{Failure conditions} We add the following failure conditions. Numerical simulation stop if any of those conditions is satisfied:
	\begin{itemize}
		\item Maximal solving time: runtime over 1.5 hour.
		\item Maximal number of iterations: failed and successful steps $ > I_{max}=10^5$.
		\item Maximal number of failed steps: number of failed steps $> 0.85\;I_{max}=8.5\times 10^{4}$.
		\item Slow solver: runtime over 45 min and less than $10\%$ of the simulation completed.
		\item Step size too small: To prevent the solver from reaching the magnitude limit of the finite-precision format we add the following condition $h>h_{min}=100\epsilon_m$ where $\epsilon_m$ is the machine epsilon (the smallest positive value summable to 1).
	\end{itemize}

	\subsection{Mixed-precision for stage computations}\label{SbS:MP-RK}
	We introduce the mixed-pre\-cision into the stage computations (\Cref{Tab:MP-Distrib}). 
	For this study we only mix single (S) and double (D) arithmetic precision formats from the IEEE 754 standard \cite{ZURAS-2008}.
	We make the following choices:
	\begin{itemize}
		\item All the operations coming from the linear combinations of the numerical stages are performed in high precision (D).
		\item For all mixed-precision solvers, the solution is stored in high precision (D).
		\item For a given stage we can choose different precision for the evaluation of $G_{ij}$, $F_{i}$ and the average of the interactions ($G_{ij}$), \ie{} three arithmetic precision can be chosen for one stage.
		\item In the case of accumulation in high precision (D) of low precision (S) terms, low precision terms are cast into high precision before accumulation.
		\item We used safeguard casts to ensure the adequate precision of the computations. Consequently, the real performance does not represent the theoretical one.
	\end{itemize}
	We use two versions of mixed-precision distribution inside the numerical scheme (denoted Mixed1 and Mixed2).
	In addition to Mixed1 and Mixed2, we also compute solutions with mono-precision solvers: single precision (Single), and double precision (Double).
	\begin{table}[htbp]
		\begin{center}
			\begin{tabular}{c|ccc|ccc|ccc|ccc}\hline
				&\multicolumn{3}{c|}{Single}&\multicolumn{3}{c|}{Mixed1}&\multicolumn{3}{c|}{Mixed2}&\multicolumn{3}{c}{Double}\\
				Stage& $ F $ & $\sum$&$G$ &$ F $ & $\sum$&$G$ & $ F $ & $\sum$&$G$ & $ F$ & $\sum$&$G$\\ \hline\hline
				$k_2$& S & S & S & S & S & S & D & D & S & D & D & D\\
				$k_3$& S & S & S & S & S & S & D & D & S & D & D & D\\
				$k_4$& S & S & S & D & D & S & D & D & S & D & D & D\\
				\hline
			\end{tabular}
			\caption{Mixed-precision distribution among the Bogacki-Shampine 3(2) pair. F corresponds to the precision used for the evaluation of $F_i$, $\sum$ for the accumulation of the interaction terms and $G$ for the interaction $G_{ij}$. S stands for single precision and D for double precision.}		
			\label{Tab:MP-Distrib}
		\end{center}
	\end{table}

	\subsection{Error metric}\label{SbS:Error-Metric}
	We present three metrics used for assessing the accuracy of the four solvers. These metrics are defined for tests that are terminated by all solvers, \ie{} the solution covers the whole integration interval. If one solver did not finish a test due to any reason among the stop criteria (see paragraph "Failure conditions" in \Cref{SbS:MP-Numerical-Scheme}), then this test is excluded from analysis for all solvers.
	
	\paragraph{Normalized final error}
	A reference solution, $X_{ref}$, is computed in double precision with the Matlab solver \texttt{ode45} which uses the 5(4) Dormand-Prince pair \cite{SHAMPINE-1997} with a $10^{-9}$ relative tolerance. The normalized final error is
	\begin{equation}\label{EQ:Metric}
		|| X_{ref}(t_{f})-X(t_{f})||_{X} = \frac{|| X_{ref}(t_{f})-X(t_{f})||_{2}}{\sqrt{N}},
	\end{equation}
	where $t_{f}$ is the final time in the test and $N$ the number of agents (see \Cref{EQ:General-Model}). The normalization by the size of the system enables to facilitate the comparison of the error on each element with respect to the theoretical error and for different sizes.
	
	\paragraph{Real local error} For \LCO{} (\Cref{EQ:LCO}), it is possible to compute an analytic solution and evaluate the real error. The real local error is the difference between analytic solution and the numerical solution after one time step.
	
	\begin{equation}\label{EQ:LCO-Local-Error}
		E_{analytic}^n = \max_{k}\Big(\frac{|X_{k}^{n+1}-X_{k,ex}(X^n,h_n)|}{max\big(|X_{k,ex}(X^n,h_n)|,\delta\big)}\Big), \quad \mathrm{with} \; h_n = t_{n+1}-t_n, \, k \in \{1,..., dN\},
	\end{equation}
	where $X_{ex}$ is the analytic solution at time $h_n$ with $X_n$ as initial condition and computed in high precision. As the initial condition for the analytic solution is always updated to the last step solution, the error accumulated in previous steps have no impact on the calculation. We took a safeguard $\delta = \frac{Abs}{Rel}$, ratio of tolerances in order to stick to the local error approximation computed by the solver, see \Cref{EQ:ERROR-O23}.
	
	\subsection{Performance metric}\label{SbS:Perf-Proxy}
	The real performance of a mixed-precision method is based on the expected computational speed-up per computation unit (function evaluation, or arithmetic operations), times the number of units to be computed. The speed-up associated with low-precision includes memory movement and arithmetic operations which are the two main sources of numerical performance bottleneck, and is implementation and architecture-dependent. When using adaptive schemes, the number of units to compute depends on the solution accuracy through the number of steps, and ultimately on arithmetic precision. To reflect the constraints we introduce a performance metric that depends on an expected speed-up coefficient $r$, and the total number of numerical steps used to solve the problem.

	Let $t_h$ and $t_l$ be the times, in high and in low precision, covering the costs for one floating-point operation as addition or multiplication, including the memory movement of the data, reading, operation itself and results writing. We define the time ratio with respect to the highest precision ($t_h$), used as reference:
	\begin{align}\label{EQ:Ratio-Proxy}
		r:=\frac{t_l}{t_h}.
	\end{align} 
    We assume that there is an expected speed-up for floating-point operation when using low arithmetic precision instead of high one \cite{BABOULIN-2009}. This speed-up may come from several hardware aspects. Smaller formats are faster to compute, to access and move in memory. In addition, optimizations such as vectorization or SIMD (Single Instruction, Multiple Data) are more efficient due to the smaller memory size of lower precision formats. This leads to:
	\begin{align*}
		0<r<1.
	\end{align*}
	We denote the total computing times $T_h$ and $T_m$, in high precision and in mixed-precision respectively. We also defined their ratio as follows:
	\begin{align}
		\Gamma := \frac{T_m}{T_h}.
	\end{align}
	
	As a measure of performance we look at the number of floating-point operations. We introduce $\Theta$, the total number of operations realised in one solver step. It is composed of the evaluation of the ODE system function, $\Theta_{ev}$ and the numerical scheme combinations $\Theta_{sch}$.
	\begin{align*}
		\Theta_{ev} &= \theta_F(d) N +\theta_G(d) N^2 + \theta_{w}(d)N^2,\\
		\Theta_{sch} &= \theta_{sch}(s)\;dN,\\
		\Theta &= s\;\Theta_{ev}+\Theta_{sch}=s N^2(\theta_G(d) + \theta_{w}(d))+ (s\theta_F(d) +d\theta_{sch}(s))N,
	\end{align*}
	where
	\begin{itemize}
		\item $\theta_F(d)$ and $\theta_G(d)$ are the numbers of floating-point operations required for evaluating respectively the functions $F_i$ and $G_{ij}$ (Section \Cref{SbS:ODE-Systems}). We emphasise that the number of operations is impacted here by the dimension of the agent $d$ and not of the size of the population $N$. 
		\item $\theta_{w}(d)$ is the number of floating-point operations for the weighting of the interactions (Hadamard product). The factor $N^2$ accounts for the summation of all the interactions for every agent.
		\item $s$ is the number of stages required by the numerical scheme and $\theta_{sch}(s)$ the number of linear combinations between the stages inside the numerical scheme (for a scalar problem).
		\item It can be seen that for large $N$, $\Theta \sim sN^2$. We assume that the majority of the computational cost lies in the stage evaluation, especially in the interaction terms.
	\end{itemize}
	
	For the whole computation this has to be multiplied by the number of steps required. We introduce the mixed-precision influence with 3 parameters: $\varrho$, $\beta$ and $\gamma$. The first one is the fraction of operations performed in lower precision. The second, $\beta$, is the ratio between the number of steps (successful and failed are included) computed in high precision and in mixed-precision, so it takes into account impacts on the global integration trajectories indicated by the successful step number and also the efficiency of the adaptive process through the failed step number\footnote{If $\beta$ is far from 1, one will not directly see if it is due to a very different integration trajectory or to a larger number of failed steps.}. The last one, $\gamma$ is the ratio of computing times for one step between mixed-precision solver and high precision solver which include the memory movements and the mathematical operations. We can express the ratio $\Gamma$ between the computing times in mixed-precision and high precision as the ratio between $\gamma$ and $\beta$:
	\begin{align}
		T_h  &= N_{step,h} \Theta \; t_h, \nonumber \\
		T_m  &= N_{step,m} \Theta (\varrho t_l + (1-\varrho) t_h), \nonumber \\
		\beta &= \frac{N_{step,h}}{N_{step,m}}, \label{EQ:Beta-2MP} \\
		\gamma &= \frac{\varrho t_l + (1-\varrho) t_h}{t_h}=\varrho r +(1-\varrho),\label{EQ:Gamma-2MP} \\
		\Gamma &= \frac{T_m}{T_h} = \frac{\gamma}{\beta}.\label{EQ:GammaMax-2MP} 
	\end{align}
	A ratio $\Gamma<1$ means that the mixed-precision solver is faster than the high precision solver.
	
	\subsection{Descriptive statistics}
	
	For each representation of the results, we select only tests that have been finished by all the solvers (\SPS{}, \MPA{}, \MPB{}, \DPS{} and the reference solver). We use some statistical quantities:
	
	\begin{itemize}
		\item The median as a reference value with 5-th and the 95-th percentiles as error bars for representing the evolution of the normalized final error against other variables (size and relative tolerance).
		\item The first, the second and the third quartiles for the box plots with the 1-st and the 99-th percentiles as whiskers.
		\item For figures with local error, we have values for each time step for every selected test. We take the averaged value of the local errors ($E_{analytic}^n$ and $Err_{BS}^n$). In that way, we have one value for each variable by test. In the figures we plot the median, with 5-th and 95-th percentiles as error bars for the averaged variables (approximated and real local error).
	\end{itemize}

\section{Results and discussion}\label{S:Results}
	All the code is available in a public repository on Inria's GitLab\footnote{Link to Inria's Gitlab repository: https://gitlab.inria.fr/amarzora/mixed-precision-adaptive-solver.}.
	\subsection{Application of the performance metric}\label{SbS:Perf-Criteria}
	 First, we look at the distribution of $\beta$ (\Cref{EQ:Beta-2MP}) among the finished tests \figcom{see \Cref{fig:NBTEST}  \Cref{AS:Tests-selection}, for the number of failed tests}{}, \ie{} the ratio of number of steps between \DPS{} and the mixed-precision solvers.
	
	\begin{table}[htbp]
		\centering
		\resizebox{\columnwidth}{!}{
		\begin{tabular}{c|c|c|c}
			& \LCO{} & \Kur{} & \CC{}\\
			\hline
			Single& (0.98, 1.0, 1.0, 1.0, 1.026)& (0.427, 0.997, 1.0, 1.0, 1.153)& (0.979, 1.0, 1.0, 1.0, 1.017)\\
			Mixed1& (0.987, 1.0, 1.0, 1.0, 1.015)& (0.482, 0.997, 1.0, 1.0, 1.118)& (0.986, 1.0, 1.0, 1.0, 1.017)\\
			Mixed2& (0.989, 1.0, 1.0, 1.0, 1.016)& (0.507, 0.998, 1.0, 1.0, 1.103)& (0.996, 1.0, 1.0, 1.0, 1.003)\\
		\end{tabular}
		}
		\caption{Distribution of $\beta$ parameter for each solver by benchmark. Each vector is composed as (min, 1st quartile, median, 3rd quartile, max) of $\beta$ values for one solver under one benchmark. The statistics are done over the tests that have been finished by all the solvers.}\label{TAB:Beta-values}
	\end{table}
	
	\Cref{TAB:Beta-values} shows indicative values (minimum, 1st quartile, median, 3rd quartile and maximum) of the $\beta$ parameter. The first remark is that for all benchmarks, all solvers perform as many steps as the \DPS{} solver in the majority of the tests. In a very few tests and only for \Kur{} the ratio is far from 1. Indeed, for all the benchmarks the 1st and third quartiles are very close or equal to 1.
	
	As the number of steps for an adaptive scheme is not fixed before the solving, we want to estimate if the difference between a mixed-precision and double precision solver in terms of supplementary steps is not large enough, $\beta \ll1$, to balance the computing speed-up of low precision ($\gamma<1$). In this situation, the ratio of the total computing time between mixed and high precision, $\Gamma$, is maintained under 1, \ie{} a better performance for the mixed-precision solver. That is expressed by the following relation between the parameters:
	\begin{align*}
		\Gamma = \frac{\gamma}{\beta} < 1  \Leftrightarrow \gamma < \beta.
	\end{align*}
	A better performance is observed for all tests only if the lowest $\beta$ is higher than $\gamma$. If the ratio $r$  between the computing time for floating-point operations in low and high precision (\Cref{EQ:Ratio-Proxy}) was fixed and not impacted by any internal supplementary cost such as casting of variable, we would be able to estimate the value of $\gamma$  (\Cref{EQ:Gamma-2MP}). That is why we consider a "naive" approach whereby that single precision operation is twice as fast as a double operation. For each mixed-precision solver, we estimate the fraction of low precision operations, $\varrho_k$. We choose to perform all the operations of the scheme, $\Theta_{sch}$, in high precision.
	Considering that systems are large ($N \gg 1$), we compute a limit value $\varrho_k$ as follows (for each version, refer to \Cref{Tab:MP-Distrib}): 
	\begin{align*}
		\text{Version 1 (\MPA{}):}\;& \varrho_1 = \frac{2\Theta_{ev}+N^2\theta_{G}(d)}{\Theta}, & &\lim_{N\to + \infty} \varrho_1 = 1- \frac{\theta_{w}(d)}{3(\theta_G(d) + \theta_{w}(d))}.\\
		\text{Version 2 (\MPB{}):}\;& \varrho_2 = \frac{3 N^2\theta_{G}(d)}{\Theta}, & &\lim_{N\to + \infty} \varrho_2 = 1- \frac{\theta_{w}(d)}{\theta_G(d) + \theta_{w}(d)}. 
	\end{align*}
	Finally, we consider a limit case  $\theta_{w}(d)\simeq\theta_{G}(d)$. The complexity of $G_{ij}$ increases the value of $\theta_G(d)$ whereas the weighting of the interactions increases only with the dimension of the agent. Normally, in models of interest we would have $\theta_G(d) \leq \theta_w(d)$, leading to a higher  $\varrho_k$. Yet, with this consideration we obtain the following values of $\varrho_k$ and $\gamma_k$ respectively for \MPA{} and \MPB{}. The trivial case for \SPS{} is also indicated:
	\begin{align*}
		\text{Version 0 (\SPS{}):} \quad \varrho_0 = 1 \implies& \; \gamma_{0} = r = 0.5,\\
		\text{Version 1 (\MPA{}):} \quad \varrho_1 = \tfrac{5}{6} \implies& \; \gamma_{1} = \tfrac{7}{12}\simeq 0.583,\\
		\text{Version 2 (\MPB{}):} \quad \varrho_2 = \tfrac{1}{2} \implies& \; \gamma_{2} = \tfrac{3}{4} =0.75.
	\end{align*} 
	$\gamma_{2}$ is higher than $\gamma_{1}\simeq 0.583$ due to fewer operations performed in low precision.
	
	For \LCO{} and \CC{}, the minimal values for $\beta$ are very close to 1, around 0.98, for all solvers. That means in the worst case, our mixed-precision solvers perform at most $2\%$ more steps than the \DPS{}. But as we stated before, this will not be enough to balance the gain $\gamma_{k}$ that we estimate reachable with low precision in computational time through the $r$ parameter.
	
	For \Kur{}, the minimal values are very low, and under the critical value of $\gamma_{k}$ for all solvers. 0.43 for \SPS{}, 0.48 for \MPA{}, and 0.51 for \MPB{}. To keep $\Gamma \leq 1$ in these cases, the parameter $r$ should be smaller than $0.5$ (below $ 0.42$ for \SPS{}). That is more restrictive but not impossible as the single precision takes up half as much memory space as double precision \cite{BABOULIN-2009}. In addition, for the last two benchmarks, the $G_{ij}$ functions are not linear. Thus, the complexity of $G_{ij}$ functions leads to a larger $\theta_{G}$, increasing the ratio of low arithmetic operations $\varrho_k$ and lowering the time ratio for one step $\gamma_{k}$.
	
	However, one should pay attention to the benchmark and the numerical scheme which are used. Indeed, our approach can be adapted to any numerical scheme, but the results will depend on the solution path taken by the solver. This can be impacted by the non-linearity and the complexity of the benchmarks but also by the choice of the numerical scheme. This is caught by $\beta$ parameter, which in our case seems clearly negligible.
	
	Another aspect is the dependence of the hardware. That is included in $r$ and $\gamma$ parameters. The $r$ parameter includes computational performances such as vectorization, parallelization or memory access which benefit for lower formats, but $r$ could be also negatively impacted by the cost of variable casting. We did not push the study on this aspect due to the difficulty to access how floating-point operations are performed.
	
	As we indicated in the mixed-precision implementation, the solution is stored in high precision (D), and all the computations related to the linear combination of the scheme just as the approximation of the local error are performed in this precision. We anticipate that for large sizes, the dominant cost should arise from interaction terms, that is why we inserted low precision in their evaluation. Yet, a possible modification would be to insert low precision in the other parts of the solver in order to modify the $\rho$ parameter and allow greater exploitation of the computational benefits captured by the $r$ parameter.

\subsection{Accuracy vs tolerance, for a fixed system size} \label{SbS:Accuracy-Tolerance}
In this section, we compare the evolution of the normalized final error (\Cref{EQ:Metric}) with respect to the relative tolerance imposed on the solver, at fixed system size. As the tolerance is reduced, \ie{} becomes stricter, the normalized final error is expected to decrease. It is also expected that for very tight tolerances, the machine epsilon of the single precision format (\Cref{SbS:MP-Numerical-Scheme}) limits the error reduction for \SPS{}. In this situation, we give a closer look at \MPA{} and \MPB{} errors, and study how the presence of low precision computation limits their accuracy. 
	
	\begin{figure}[htbp]
		\centering
		\includegraphics{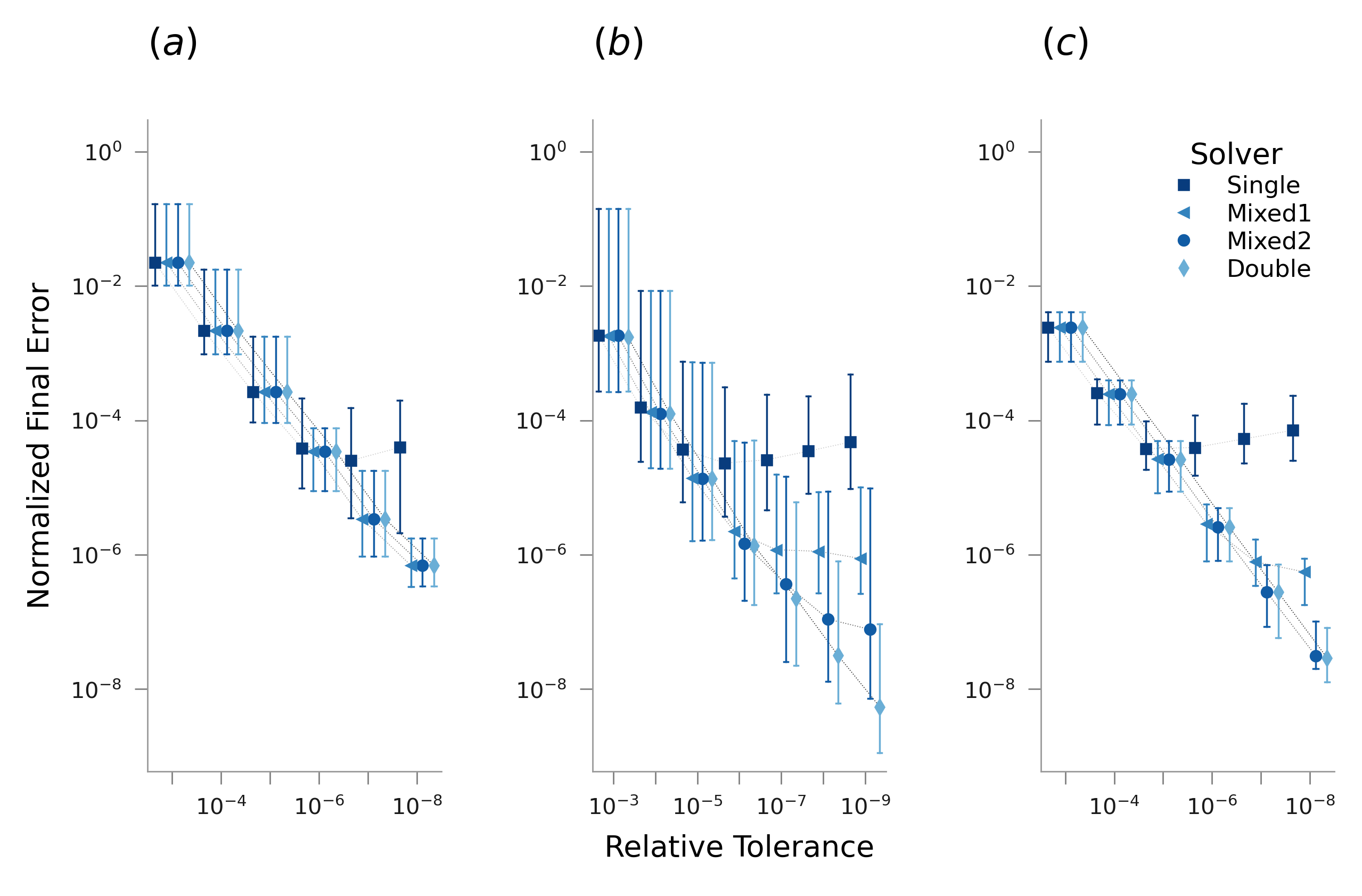}
		\caption{Normalized final error, see \Cref{EQ:Metric}, with respect to Matlab's \texttt{ode45} solution. 1000 agents for \LCO{} (a), to 2000 for \Kur{} (b) and to 700 for \CC{} (c). The error bars are the 5-th and 95-th percentiles. For tolerances of ($10^{-3}$, $10^{-4}$, $10^{-5}$, $10^{-6}$, $10^{-7}$, $10^{-8}$, $10^{-9}$), there are respectively in (a) (113, 97, 78, 59, 33, 19, 0), in (b) (498, 423, 359, 270, 189, 122, 53) and (c) (289, 235, 199, 134, 96, 42, 0) performed tests. Each solver is identified by its symbol: \SPS{} square, \MPA{} triangle, \MPB{} circle and \DPS{} diamond. As the tolerance values tested are discrete, a small shift is added between the solvers to distinguish their symbol and error bars.}\label{fig:Tol-VS-Error}
	\end{figure}
	
For the 3 benchmarks and all solvers, the normalized final error decreases as the tolerance is reduced but stays at large values ($10^{-3}$, $10^{-4}$, $10^{-5}$). 
As expected, the normalized final error of \SPS{} \figcom{\Cref{fig:Tol-VS-Error}}{squares} plateaus for relative tolerances below $10^{-5}$ in all benchmarks. Indeed, the normalized final error is between $10^{-4}$ and $10^{-5}$ for the three benchmarks. Even a slight rise is visible while the relative tolerance continues to be reduced. The accuracy is limited by the machine epsilon in single precision ($10^{-8}$), \SPS{} solution does not satisfy the tolerance criterion in the sense that there is no difference in the error for stricter tolerances.

In contrast, errors steadily decrease in the three other solvers, albeit with noticeable differences. In the case of \LCO{} \figcom{\Cref{fig:Tol-VS-Error}}{a}, the averaged normalized final errors coincide for all three solvers. The mixed-precision solvers have an accuracy equivalent to \DPS{}. However, even if \MPA{} \figcom{\Cref{fig:Tol-VS-Error}}{triangles} has a better accuracy than \SPS{}, its error is also limited at stricter tolerances, below $10^{-7}$. The normalized final error plateaus around $10^{-6}$ for \Kur{} and \CC{} \figcom{\Cref{fig:Tol-VS-Error}}{b and c}, whereas \MPB{} and \DPS{} \figcom{\Cref{fig:Tol-VS-Error}}{circles and diamonds} have a normalized final error that reaches $10^{-8}$ and does not plateau. Yet, on \Kur{}, \MPB{} seems to plateau from tolerances $10^{-8}$ and smaller, its normalized final error is better than \MPA{} error, but it is approximately $10$ times larger than \DPS{} error.
	
These results highlight the limitations of low precision operations in terms of accuracy, especially when dealing with stringent tolerance requirements. While single precision algorithms offer computational efficiency, their limited potential for high accuracy under tight tolerances must be carefully considered in applications where accuracy is a strong criterion. Mixed-precision methods have the possibility to use the computational efficiency of low arithmetic precision, while pushing back the accuracy limit even for stringent tolerances.
	
\subsection{System size vs accuracy, for fixed relative tolerance}\label{SbS:Accuracy-Size}
In this section, we test how the system size affect the accuracy of the solution. We look at the evolution of the normalized final error (\Cref{EQ:Metric}) with respect to the size of the system. The relative tolerance is fixed. Because we want to use mixed-precision as a lever for modeling large-scale systems, we have to check that when the scale is increased, the low precision does not become prohibitive.
 
	\begin{figure}[htbp]
		\centering
		\vspace{-60pt}
		\includegraphics{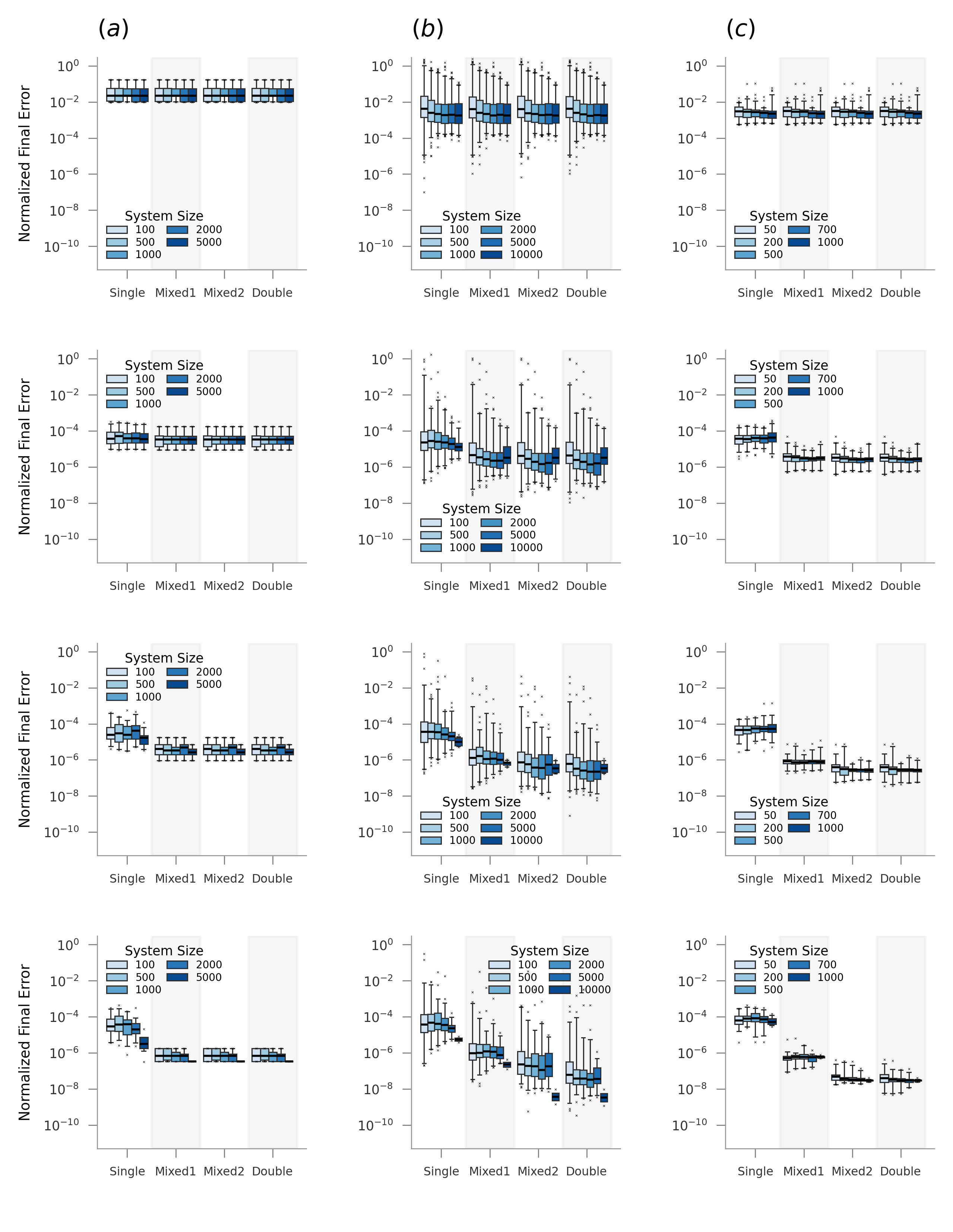}
		\caption{
			Distribution of normalized final error, see \Cref{EQ:Metric}, for all the tests completed by all the solvers. 
			One box plot represents a distribution, with limits at 1-st and 99-th percentiles, for all the tests of one benchmark, indicated by the column (a) \LCO{}, \Cref{EQ:LCO}, (b) \Kur{}, \Cref{EQ:Kuramoto} and (c) \CC{}, \Cref{EQ:CC}. With a unique system size, indicated by the intensity of blue, pale blue for small sizes (100 for \LCO{} and \Kur{} and 50 for \CC{}) to dark blue for the largest ones (5,000 for \LCO{}, 10,000 for \Kur{}, 1,000 for \CC{}). At a specific relative tolerance that can be found with the row, from top to bottom ($10^{-3}$, $10^{-6}$, $10^{-7}$, $10^{-8}$). Run by one solver, indicated on the abscissa of the boxes group, from left to right (\SPS{}, \MPA{}, \MPB{}, \DPS{}).
		}
		\label{fig:Size-VS-Error}
		
	\end{figure}
	
	The median \figcom{\Cref{fig:Size-VS-Error}}{black horizontal lines} and the first and 99th percentiles of the normalized final error are either unchanged or tend to decrease as the system size increases. For example, with \LCO{}, \figcom{\Cref{fig:Size-VS-Error}}{column a}, there is no real change as the size increases. In addition, the mixed-precision solvers have a normalized final error as accurate as the \DPS{} at every tolerance. The effect of the size is more pronounced in the case of \Kur{} \figcom{}{column b}, at every tolerance the whiskers and the median are lowered as the size increases. For \CC{} \figcom{}{column c} the decrease of the error with the increasing size is present but less pronounced. 
	
	Furthermore, for all benchmarks, \MPB{} has results closed to \DPS{} ones at every tolerance. The error obtained on larger systems by \MPB{} ends up being of the same order as the error of \DPS{} for lower sizes at same tolerance. 
	For example, for \Kur{} tests solved with a $10^{-6}$ tolerance \figcom{}{second row, column b}, the median of the normalized final error for systems of 5,000 oscillators computed with \MPB{} is lower than the \DPS{} one for tests with sizes of 500 or 1,000 oscillators at the same tolerance. 
	Similarly, for \CC{} tests solved with a $10^{-6}$ tolerance \figcom{}{second row and column c}, the normalized final error for systems of 1,000 agents solved by \MPA{} and \MPB{} are equivalent or better than the \DPS{} normalized final error for systems of 100 agents. 
	
	The previous observation is model-dependent in the sense that is more pronounced for \Kur{}, \figcom{\Cref{fig:Size-VS-Error}}{column b} than for \LCO{}, \figcom{}{column a} and \CC{} \figcom{}{column c}. However, either the errors are equivalent between mixed-precision solvers and high precision solver, or, with the error improvement on large systems, the mixed-precision solvers have an error equivalent to the high precision solver error for same constraint but on smaller systems. Even if there is a small loss in accuracy, as expected with lower precision, it is limited with mixed-precision solvers and tends to be reduced as the system size increases. In addition, for large systems the potential speed-up is more important, as underlined in \Cref{SbS:Perf-Criteria}, because the fraction of low precision operations increases as the size increases. 
	
	\subsection{Adaptive scheme: a procedure limited by the arithmetic precision} \label{SbS:Local-Error}
	In this section, we look at the efficiency of the error approximation in the Bogacki-Shampine 3(2) pair with respect to the arithmetic precision. For \LCO{}, the existence of an analytic solution allows us to compare the approximate local error ($E_{BS}$,  \Cref{EQ:ERROR-O23}) between two numerical solutions computed with the scheme, and the real local error ($E_{analytic}$, \Cref{EQ:LCO-Local-Error}) between the numerical solution and the analytic one. 
	
	\begin{figure}[htbp]
		\centering
		\includegraphics{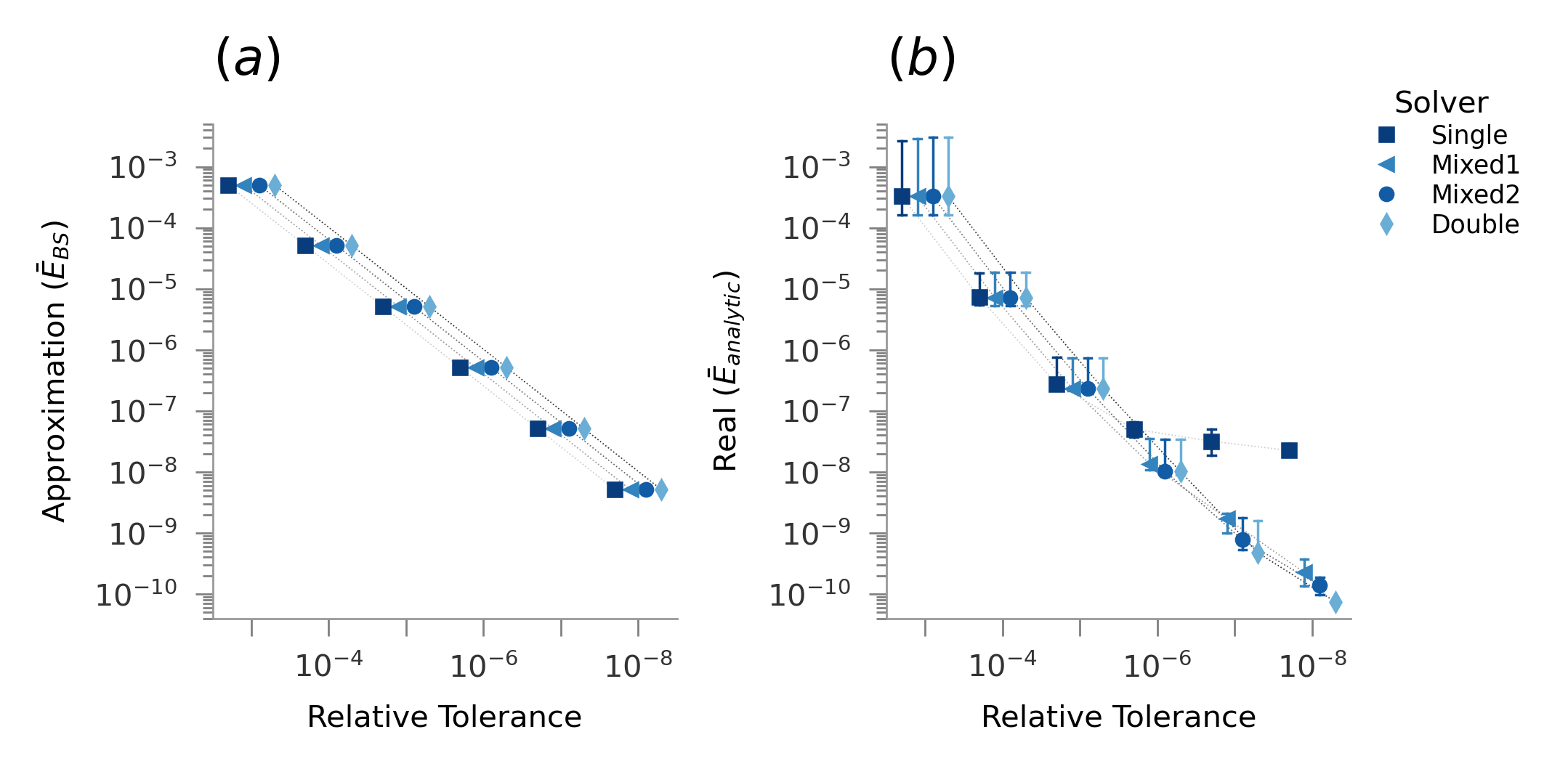}
		\caption{(a) Averaged local error (one value per test) approximated as the difference between the 2 RK schemes ($E_{BS}$, in \Cref{EQ:ERROR-O23}).
			(b) Averaged real local error (one value per test). In each plot the median which is computed over all the tests performed for a fixed tolerance is given with the error bars corresponding to the 5th and 95th percentiles. As the tolerance values tested are discrete, a shift was added to distinguish the points and the error bars. In (a) and (b) square is for \SPS{}, triangle \MPA{}, circle for \MPB{}, and diamond for \DPS{}.}
		\label{fig:LCO-Local-Error-Comp}
	\end{figure}
	
As expected, the approximated local error is reduced as the tolerance tightens \figcom{\Cref{fig:LCO-Local-Error-Comp}}{a} and its value is always lower than the relative tolerance, meaning that for all the solvers, the tolerance criteria are satisfied. However, regarding the local error computed against the analytic solution \figcom{\Cref{fig:LCO-Local-Error-Comp}}{b}, the tolerance criterion is not satisfied for \SPS{} \figcom{\Cref{fig:LCO-Local-Error-Comp}}{b, squares}. The averaged local error plateaus around $10^{-7}$ for tolerances lower than $10^{-6}$. Indeed, the machine epsilon of single precision ($10^{-8}$) limits the real deviation between the analytic solution and the numerical one. On the contrary, mixed-precision solvers real errors \figcom{\Cref{fig:LCO-Local-Error-Comp}}{b, triangles and circles} are reduced as the tolerance is tightened. The real and the approximated errors evolve similarly to \DPS{} errors \figcom{\Cref{fig:LCO-Local-Error-Comp}}{b, diamonds}, even at stricter tolerance. 

The specificity of adaptive schemes lies in the step size selection procedure. The solver computes an approximated local error in order to adapt the step size to the tolerance. In \figcom{\Cref{fig:LCO-Local-Error-Comp}}{a}, the procedure works for all solvers with respect to the approximation. However, the low precision impacts the real accuracy of the solution. In \figcom{\Cref{fig:LCO-Local-Error-Comp}}{b}, the real local error at strict tolerance is higher than the approximation of \SPS{} and the tolerance criterion is not satisfied. However, the mixed-precision solvers are not impacted in that way for the tested tolerances. Their accuracy is close to the high-precision solver.
	
\section{Conclusion}\label{S:Conclusion}
In this work, we have proposed a mixed-precision version of an adaptive RK method and studied its performance through a proxy independent of the architecture and arithmetic format used. Our approach can be generalised to any similar method such as embedded Fehlberg methods \cite{FEHLBERG-1969}, and extended to any type of numerical format (\textit{e.g.} bfloat16 or fp8). We theoretically expect a performance gain that essentially depends on the ratio of operations performed in low arithmetic precision and the difference of computational cost between high and low arithmetic precision formats.  

Contrary to a full low precision solver, we showed that mixed-precision solvers are an interesting balance between computational performance and accuracy. Global performance depends on the complexity of the systems induced by the agent interactions, and by the order of the scheme and the tolerance, which affect the number of time steps. For adaptive scheme, the number of steps cannot be set before solving, but we showed that the difference regarding the number of steps, between a standard solver in double precision and others in mixed-precision is not high enough to balance the speed-up of low precision \figcom{\Cref{TAB:Beta-values}}{}.

As the tolerance tightens, numerical convergence of RK methods implies that infinite precision solution should converge to the true solution. Yet, the low precision impacts the accuracy through round-off errors, for example \SPS{} accuracy is limited by its machine epsilon for larger tolerances \figcom{\Cref{fig:Tol-VS-Error}}{}. The real advantage is the gain of accuracy with the size. Not only mixed-precision solvers achieve finer errors at strict tolerances, but their solution accuracy is improved as the size increases \figcom{\Cref{fig:Size-VS-Error}}{}. For large systems, the accuracy of mixed-precision solvers is equivalent to the \DPS{} error for smaller systems. 
	
It is also important to emphasise the role of arithmetic precision in the adaptive procedure. Error control works well for coarse tolerances, but for more restrictive ones, the single precision solver is unable to achieve the required error magnitude. Indeed, the local error with respect to the analytic solution is higher than the fixed tolerance, \ie{} the step should be rejected, but the error approximation computed by the solver is below the tolerance \figcom{\Cref{fig:LCO-Local-Error-Comp}}{}. The mixed-precision solvers are more robust and the error procedure is efficient even at strict tolerance.

Two parameters can be adjusted to get a better performance: the ratio of low precision operations over the number of performed floating-point operations within the numerical scheme, and the performance of the architecture through the $r$ parameter that can be modified with other floating-point formats or the use of  specific hardware architectures with process units (CPU, GPU) designed for mixed-precision \cite{CHOQUETTE-NVIDIA-2020}. Here, we consider only two types of precision, high (double precision) and low (single precision), which are the most widespread and available. Nonetheless, our metric can be easily extended to any floating-point format and potentially other type of formats such as fixed-point \cite{HUBRECHT-2024}. In addition, the performance proxy can be easily adapted for several native performances. It could be interesting to first investigate the runtime performance and extend this work to other numerical methods and other finite-precision formats. Indeed, computer representation of numbers requires a balance between available range, achievable precision and performance in terms of storage and computation speed. Although the floating-point formats defined by the IEEE 754 standard are the most widely used \cite{ZURAS-2008}, others are designed to emphasise one or more of the above criterion. For example, {\it Google Brain bfloat16} has the same number of bits as the IEEE 754 half precision (16-bit format) but it has a range equivalent to that of the standard 32-bit single precision format \cite{WANG-2019}.

{\footnotesize
\setlength{\itemsep}{2pt}
\setlength{\parskip}{0pt}
\bibliographystyle{acm}
\bibliography{references}
}

\newpage
	
\appendix
\section{Parameters}\label{AS:Parameters}
	For all benchmarks we have 
	\paragraph{Tolerances} The absolute tolerance is lower or equal than the relative tolerance. 
	\paragraph{Explored parameters} All the explored parameters are generated with the \textit{sobolset} function. The value is then adapted to the parameter interval (see \Cref{Tab:Parameters-Values-LCO,Tab:Parameters-Values-Kur,Tab:Parameters-Values-CC}).
	
	\subsection{Coupled linear oscillators (\LCO{})}\label{ASbS:Parameters-LCO}
	\paragraph{Initial conditions} All the initial conditions are generated in the same way. The seed is fixed by the number identifying the test. Then the function \textit{rand} of Matlab (returns a random scalar drawn from the uniform distribution in the interval $(0,1)$) is used to create a vector of initial conditions between $(0,1)$ which is multiplied by the maximum value (here $2$).
	
		\begin{table}[htbp]
		\begin{center}
			{\scriptsize
			\begin{tabular}{c|c}\hline
				Parameters & Value(s)\\
				\hline
				Number of tests& $2000$ \\
				$N$ (number of oscillators) & $\{100,500,1000,2000,5000\}$\\
				Tolerances & $\{10^{-8},10^{-7},10^{-6},10^{-5},10^{-4},10^{-3}\}$\\
				$T_f$ (Final time) & $\{10\pi, \; 20\pi, \; 50\pi \}$\\
				Initial conditions & $[0,1]$ (uniform distribution)\\
				\hline
			\end{tabular}
			}
			\caption{Explored parameters for Coupled linear oscillators.}
			\label{Tab:Parameters-Values-LCO}
		\end{center}
	\end{table}
	
	\subsection{Kuramoto (\Kur{})}\label{ASbS:Parameters-Kur}
	\paragraph{Initial conditions} All the initial conditions are generated in the same way. The seed is fixed by the number identifying the test. Then the function \textit{rand} of Matlab (returns a random scalar drawn from the uniform distribution in the interval $(0,1)$) is used to create a vector of initial conditions between $(0,1)$ which is multiplied by the maximum value (here $2\pi$).
	
	\paragraph{Natural frequencies} The natural frequencies are generated with a seed fixed with the number of the test. The function \textit{randn} (returns a random scalar drawn from the standard normal distribution) is used to create the vector which is multiplied by the parameter $\sigma$. In other words, we have $\omega \sim \mathcal{N}(0,\,\sigma^{2})$.
	
	\begin{table}[htbp]
		\centering
		{\scriptsize
		\begin{tabular}{c|c}
			\hline
			Parameters & Value(s) \\
			\hline
			Number of tests & $12 000$ \\
			$N$ (Number of oscillators)& $\{100,500,1000,2000,5000, 10000\}$\\
			Tolerances & $\{10^{-9},10^{-8},10^{-7},10^{-6},10^{-5},10^{-4},10^{-3}\}$\\
			$T_f$ (Final time) &$\frac{4\pi}{\mathrm{med}(|W|)K+0.001}$\\
			Initial conditions & $[0,2\pi]$ (uniform distribution)\\
			$\sigma$ & $[0,1]$\\
			$\omega_i$ & $[-\sigma,\sigma]$\\
			$K$(coupling coefficient) & $[0,3\sigma]$\\
			\hline
		\end{tabular}
		}
		\caption{Explored parameters for Kuramoto.}
		\label{Tab:Parameters-Values-Kur}
	\end{table}
	
	\subsection{Circadian clock (\CC{})}\label{ASbS:Parameters-CC}
	\paragraph{Initial conditions} The seed is fixed by the number of the test. Then the four components are created in this way:
	\begin{equation*}
		(\Xag_{i,1}, \; \Xag_{i,2},  \; \Xag_{i,3},\;\Xag_{i,4})=(1,1,-1.19,-0.62)+0.2\times(p_1,p_2,p_3,p_4).
	\end{equation*}
	Where $p_k$ follows a uniform distribution between $[-0.5,0.5]$.

\begin{table}[htbp]
	\centering
	{\scriptsize
	\begin{tabular}{c|c}
		\hline
		Parameters & Value(s) \\
		\hline
		Number of tests & $5000$ \\
		$N$ (Number of cells)& $\{50,200,500,700,1000\}$\\
		Tolerances &$\{10^{-8},10^{-7},10^{-6},10^{-5},10^{-4},10^{-3}\}$\\
		$T_f$ (Final time) &$[2,3] \times T_{cycle}=[48,72]$\\
		$K$(coupling coefficient)&$\{0.001,0.1,1,10\}$\\
		$I_0$&$\{0.228249, 1.5,10\}$\\
		\hline
	\end{tabular}
	}
	\caption{Explored parameters for \CC{}}
	\label{Tab:Parameters-Values-CC}
\end{table}
	
\section{Tests selections}\label{AS:Tests-selection}
	First we present the success rate of the 4 solvers in \Cref{TAB:Rate-Succ}. The rate is lower for \Kur{} but we performed more tests with higher sizes and more strict tolerances, see \Cref{Tab:Parameters-Values-Kur}. We also note that the success rate is very close among the solvers. That is why we restrict the study of the results only on the tests that have been finished, \ie{} the final time has been reached before any failure criterion, see \textbf{Failure conditions} in \Cref{SbS:MP-Numerical-Scheme}.
		
\begin{figure}[htbp]
		\centering
		\includegraphics{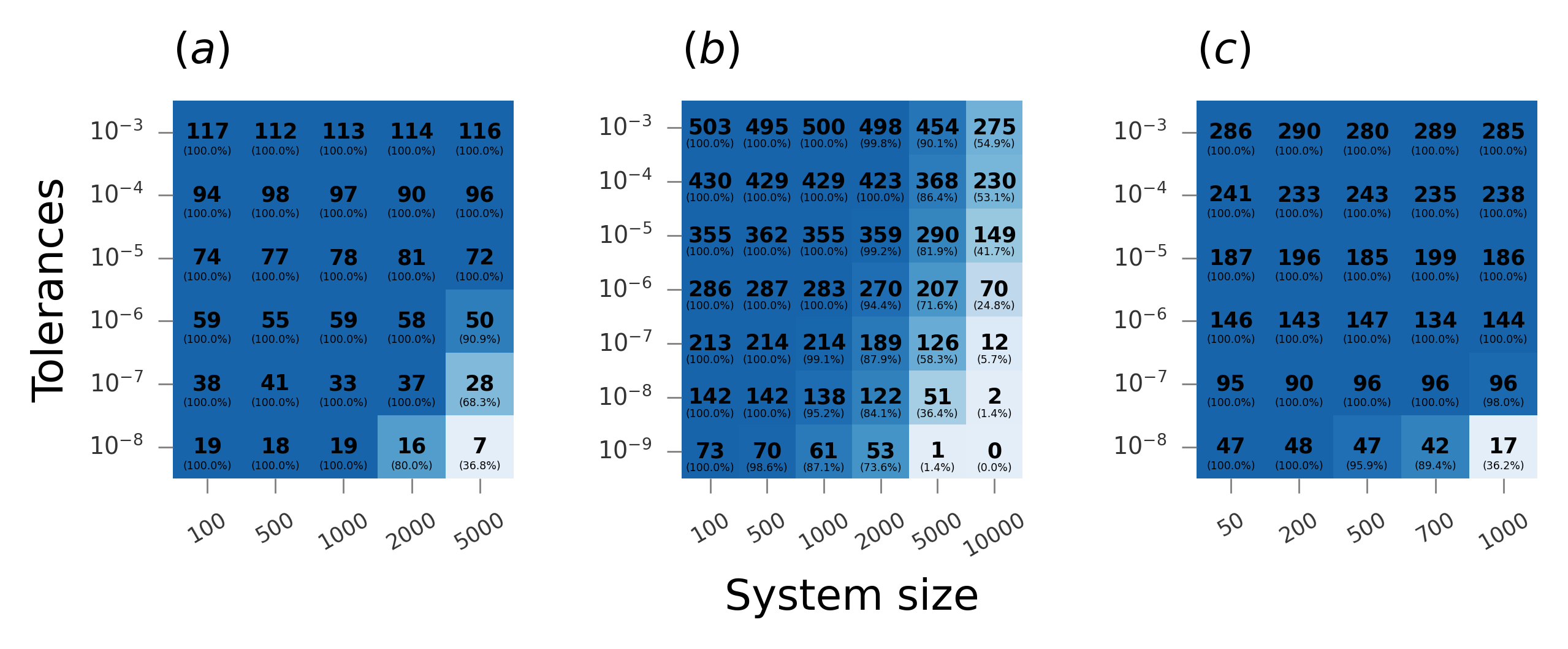}
		\caption{Number of tests successfully solved by all solvers for each couple of system size (abscissa) and relative tolerance (ordinate). (a) \LCO{}, (b) \Kur{} and (c) \CC{}.}
		\label{fig:NBTEST}
	\end{figure}

	In \Cref{fig:NBTEST} the ratio between the number of tests finished by all the solvers compared to the number of launched tests with the same couple of system size (abscissa) and relative tolerance (ordinate) and this for each benchmark.
	\begin{table}[htbp]
	\centering
	{\scriptsize
	\begin{tabular}{c|c|c|c}	
		Solver & \LCO{} & \Kur{} & \CC{} \\
		\hline
		NoT & 2,000 & 12,000 & 5,000\\
		\hline
		\SPS{} & $99.0\%$ & $87.2\%$ & $99.3\%$ \\
		\MPA{} & $98.5\%$ & $87.2\%$ & $99.3\%$ \\
		\MPB{} & $98.6\%$ & $87.2\%$ & $99.3\%$ \\
		\DPS{} & $98.8\%$ & $89.8\%$ & $99.2\%$ \\
		Ref & $98.8\%$ & $89.8\%$ & $99.2\%$ \\
	\end{tabular}
	}
	\caption{Success rate for each solver by benchmark.}\label{TAB:Rate-Succ}
	\end{table}

\end{document}